\documentclass{article}

\usepackage{arxiv}

\usepackage[utf8]{inputenc} 
\usepackage[T1]{fontenc}    
\usepackage{hyperref}       
\usepackage{url}            
\usepackage{booktabs}       
\usepackage{amsfonts}       
\usepackage{nicefrac}       
\usepackage{lipsum}
\usepackage{graphicx}
\usepackage{amsmath}
\graphicspath{ {./images/} }

\title{Frequency Distribution of Prime Numbers between an Integer and its Square: A Case Study}

\author{
 Tashreef Muhammad \\
  Department of Computer Science and Engineering\\
  Ahsanullah University of Science and Technology\\
  Dhaka, Bangladesh \\
  \texttt{tashreef.cse@aust.edu} \\
   \And
 G. M. Shahariar \\
  Department of Computer Science and Engineering\\
  Ahsanullah University of Science and Technology\\
  Dhaka, Bangladesh \\
  \texttt{shahariar\_shibli.cse@aust.edu} \\
  \And
 Tahsin Aziz \\
  Department of Computer Science and Engineering\\
  Ahsanullah University of Science and Technology\\
  Dhaka, Bangladesh \\
  \texttt{tahsinaziz.cse@aust.edu} \\
  \And
 Mohammad Shafiul Alam \\
  Department of Computer Science and Engineering\\
  Ahsanullah University of Science and Technology\\
  Dhaka, Bangladesh \\
  \texttt{shafiul.cse@aust.edu} \\  
}

\begin{document}
\maketitle
\begin{abstract}
The chronicle of prime numbers travel back thousands of years in human history. Not only the traits of prime numbers have surprised people, but also all those endeavors made for ages to find a pattern in the appearance of prime numbers has been captivating them. Until recently, it was firmly believed that prime numbers do not maintain any pattern of occurrence among themselves. This statement is conferred not to be completely true. This paper is also an attempt to discover a pattern in the occurrence of prime numbers. This work intends to introduce some mathematical well-known equations that point to the existence of a simplistic pattern in the number of primes within the range of a number and its square. We assume that the rigorous evaluation of the perceived pattern may benefit in many aspects such as applications of encryption, algorithms concerning prime numbers, and many more.
\end{abstract}


\section{Introduction}
The ancient Greek mathematicians were the first ones to inquire broadly about prime numbers and their characteristics.  The mathematicians of Pythagoras's school (500 BC to 300 BC) conceived the idea of a number being prime and enticed in \textit{perfect} as well as \textit{amicable} numbers \cite{jjo'connorefrobertson2018}. Mathematicians throughout the times have adjudicated vigorously to unhitch mysteries about prime numbers. As a result, some properties of prime numbers like twin primes and cousin primes are distinguished. However, the properties are, the concept of pattern in the number of primes has always perplexed many. 

\begin{figure}[ht] 
    \centering
    \includegraphics[width=\textwidth,height=\textheight,keepaspectratio]{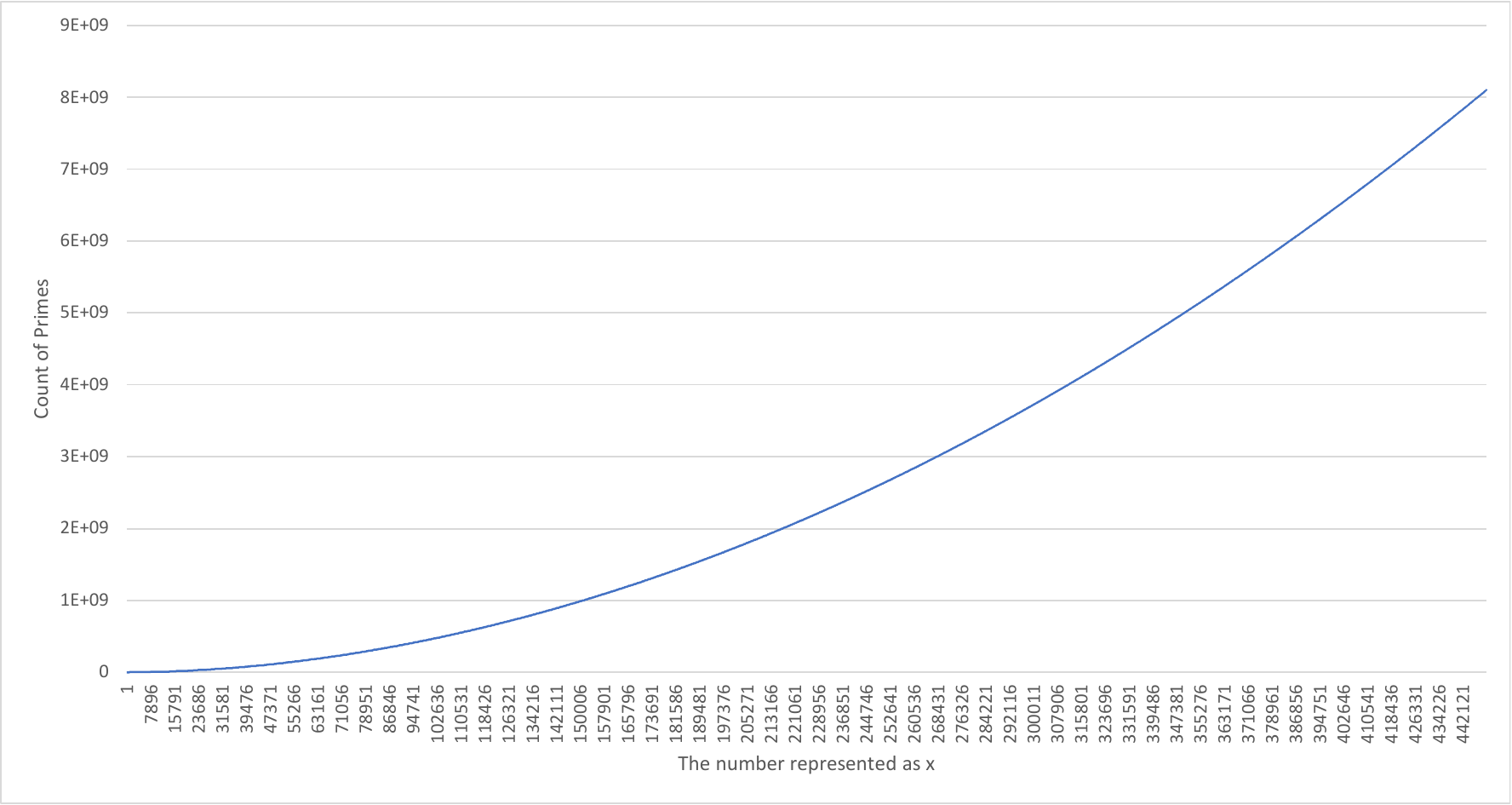}
    \caption{\textbf{Count of Primes Vs. Value of x}}
    \label{fig:cnt_prime_vs_num}
\end{figure}

Prime numbers were discovered a long time ago, but it is only in recent years that a pattern is observed in the last digit of prime numbers. In \cite{Lemke}, it is affirmed that the last digit of a prime number tends to avoid recurring in consecutive primes. Many people refer it as "conspiracy among primes" \cite{Klarreich}. This paper explains the observation that the number of times prime numbers befalling in a specific range can follow a pattern. The following claim can be observed in Figure \ref{fig:cnt_prime_vs_num}.

\begin{center}
In the range $[x, x^2]$; where $x \in  \mathbf{N}$,\\  
the number of primes occurring follows a pattern that can be represented by a continuous curve.  
\end{center}

\section{Contribution}
\begin{figure}[ht]
    \centering
    \includegraphics[width=\textwidth,height=\textheight,keepaspectratio]{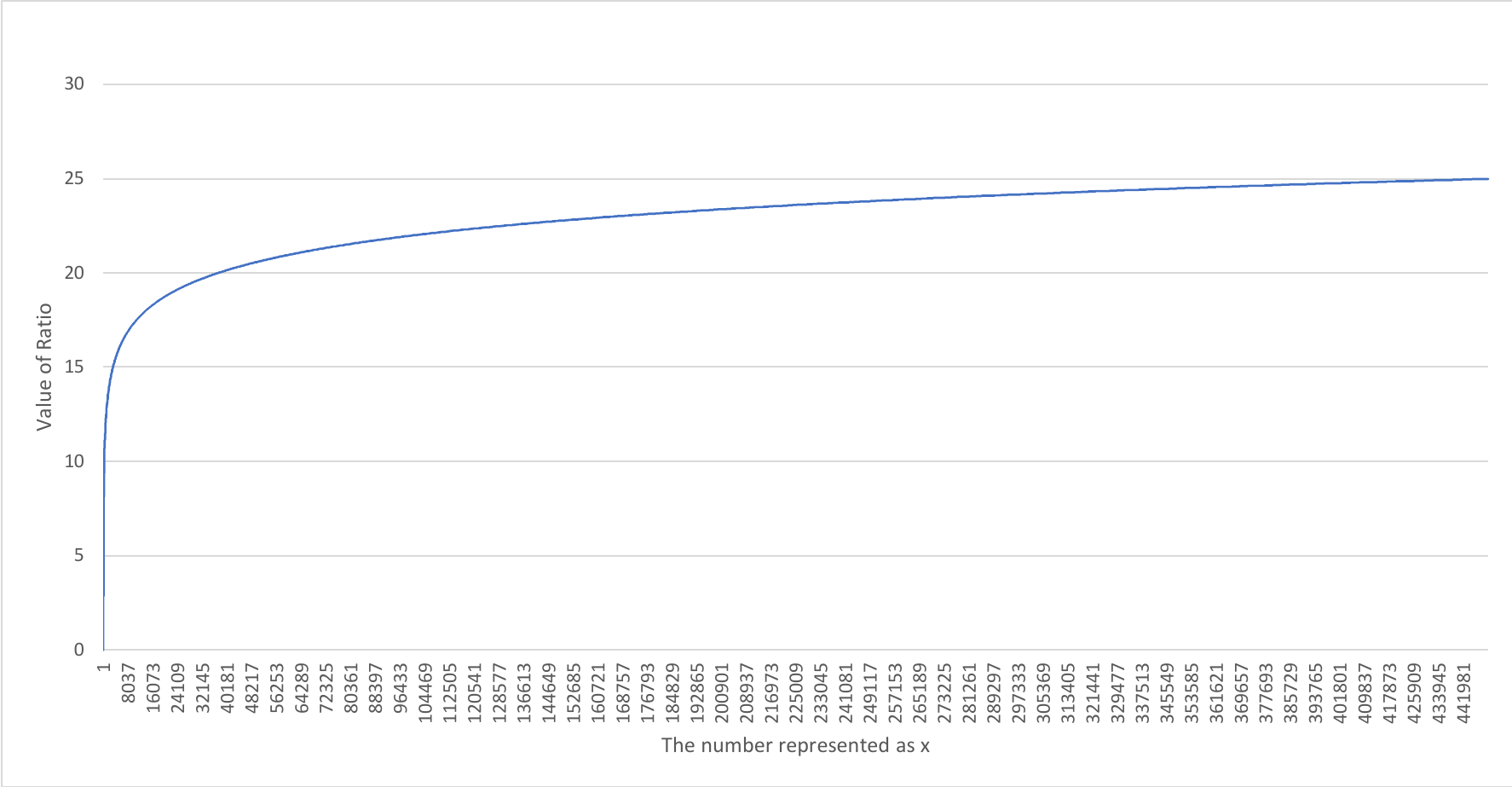}
    \caption{\textbf{Count Ratio Vs. Value of x}}
    \label{fig:ratio}
\end{figure}

The curve depicted in Figure \ref{fig:cnt_prime_vs_num} turns out to be perpetual which implies that there is a higher probability of denoting it using an equation. The observation after further research is that the ratio of range interval that is, $x^2 - x$, and the number of primes within that range also appears as a continuous curve. This curve represented in  Figure \ref{fig:ratio} resembles almost similar to a natural logarithmic curve. These remarkable details determine the presence of a somewhat pattern obvious in prime numbers.

Due to computational limitation, the pattern is only witnessed from $x = 2$ to $x = 449999$. $x = 1$ is intentionally left out in some cases as no primes fall in that range and produce undefined values. Nevertheless, bypassing only one value has a minor impact on the overall measurement. The outcome of the Figure \ref{fig:cnt_prime_vs_num} and Figure \ref{fig:ratio} guides to a foresight that presumably the pattern will pursue higher values. 

The most alike mathematical analysis related to this subject is to predict a lower limit of the number of primes in the range $[x, x^2]$ derived from Bertrand's Postulate \cite{reserachgate}. However, the derived equation provides a very low estimation. The observation is that the larger the numbers grow, the equation generates the higher Relative Error, and soon reaches $100\%$.

\begin{figure}[ht]
    \centering
    \includegraphics[width=\textwidth,height=\textheight,keepaspectratio]{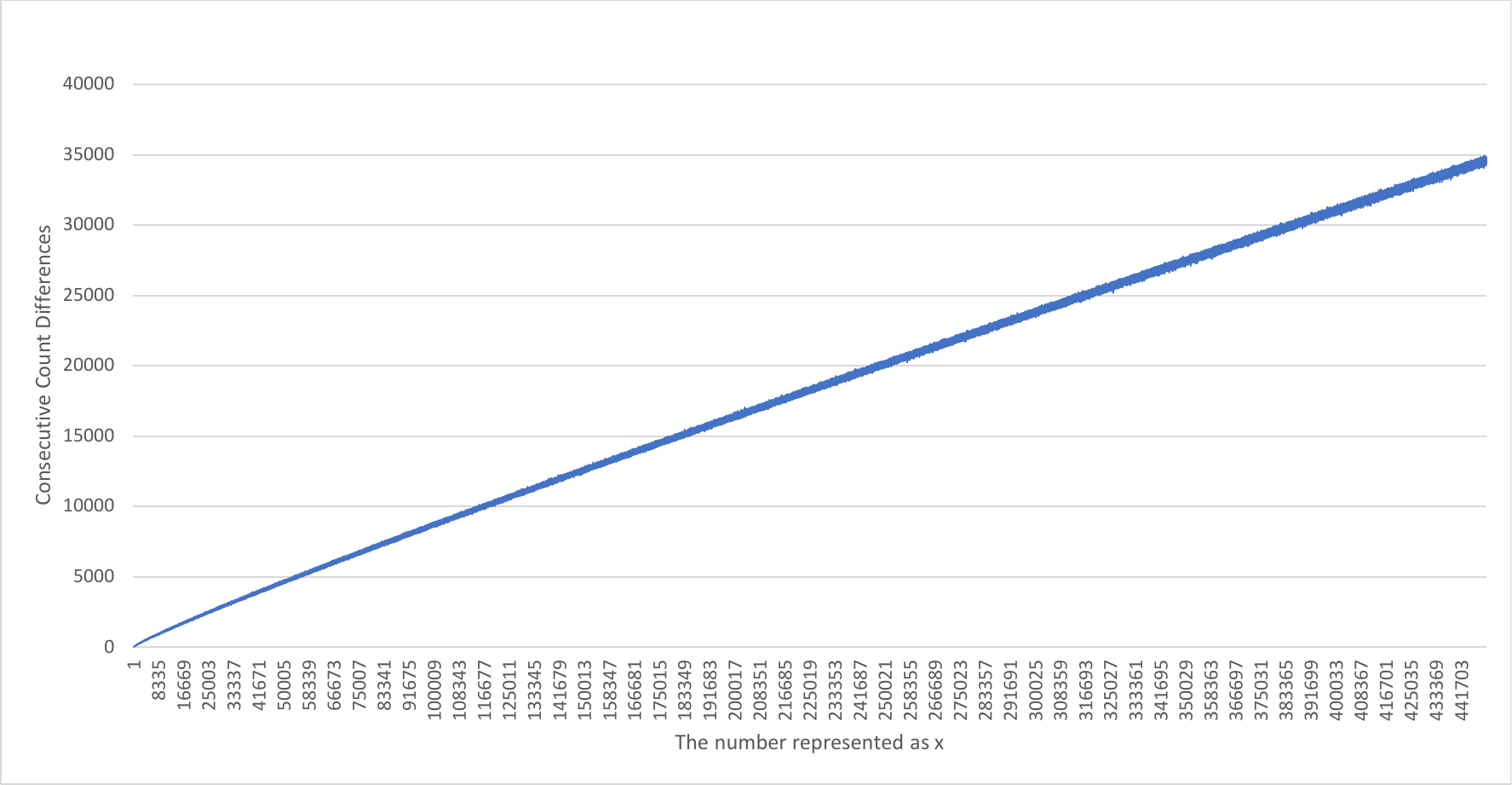}
    \caption{\textbf{Count Difference Vs. Value of x}}
    \label{fig:diff}
\end{figure}

The diversity of prime count between two contiguous ranges exposed a property that is not so arbitrary. The deviation of values improved as the number of entries gets higher. However, the line in Figure \ref{fig:diff} is approximately linear and follows a clear upward trend.

This paper also presents some approaches to predict the equation that satisfies the graphs in Figure \ref{fig:cnt_prime_vs_num}, Figure \ref{fig:ratio} and Figure \ref{fig:diff}. The ratio curve visibly represents $logarithmic$ characteristics. But for the prime count curve, multiple well-known functions are tried to fit for prediction such as $Power$ $Series$ function, $Hyperbolic$ function, $Second$ $Degree$ $Polynomial$ function and $Conics$ function. 

A $Custom$ $Equation$ (Equation: \ref{eqn:cus_calc}) is derived from Figure \ref{fig:ratio} to fit the number of primes count curve. Among all the techniques, the $Custom$ $Equation$ outperforms all the others with least \textbf{Average Relative Error} on all the points from  $x = 2$ to $x = 449999$ being $0.01\%$. All the other equations have crossed the main curve except the curve generated from the $Second$ $Degree$ $Polynomial$ function. It has highest \textbf{Average Relative Error} among all being $22.31\%$. The curve of $Second$ $Degree$ $Polynomial$ function always stays below the main curve that provides a much better lower limit of the main prime count curve than the derived equation from Bertrand's Postulate \cite{reserachgate}.

\section{Data Generation}
\textbf{Anaconda packages} \cite{anaconda} were used for the programming language \textbf{Python} to generate the count of primes. Different python packages including \textit{Pandas} \cite{reback2020pandas}, \textit{NumPy} \cite{harris2020array}, \textit{SymPy} \cite{10.7717/peerj-cs.103} were also used for assistance in the generation procedure. A few generated data are presented in Table \ref{tab:samplemain}. Figure \ref{fig:cnt_prime_vs_num} is obtained from the fully generated data.

\begin{table}[ht]
    \centering
    \caption{\textbf{Some Generated Prime Counts}}
    \label{tab:samplemain}
    \resizebox{\linewidth}{!}{\begin{tabular}{|c | c | c||c | c | c|}
         \hline
         \textbf{The Number} & \textbf{Square of Number} & \textbf{Primes in the Range} & \textbf{The Number} &  \textbf{Square of Number} & \textbf{Primes in the Range}\\ [0.5ex] 
         \hline\hline
         2 & 4 & 2 & 312402 & 97595009604 & 4023029104\\ 
        \hline
        3 & 9 & 3 & 312403 & 97595634409 & 4023053647\\ 
        \hline
        4 & 16 & 4 & 312404 & 97596259216 & 4023078351\\ 
        \hline
        5 & 25 & 7 & 312405 & 97596884025 & 4023102967\\ 
        \hline
        6 & 36 & 8 & 312406 & 97597508836 & 4023127709\\ 
        \hline
        7 & 49 & 12 & 312407 & 97598133649 & 4023152334\\ 
        \hline
        8 & 64 & 14 & 312408 & 97598758464 & 4023177020\\ 
        \hline
        9 & 81 & 18 & 312409 & 97599383281 & 4023201611\\ 
        \hline
        10 & 100 & 21 & 312410 & 97600008100 & 4023226447\\ 
        \hline
        11 & 121 & 26 & 312411 & 97600632921 & 4023251108\\ 
        \hline
        12 & 144 & 29 & 312412 & 97601257744 & 4023275787\\ 
        \hline
        13 & 169 & 34 & 312413 & 97601882569 & 4023300478\\ 
        \hline
        14 & 196 & 38 & 312414 & 97602507396 & 4023325121\\ 
        \hline
        15 & 225 & 42 & 312415 & 97603132225 & 4023349935\\ 
        \hline
        16 & 256 & 48 & 312416 & 97603757056 & 4023374564\\ 
        \hline
        17 & 289 & 55 & 312417 & 97604381889 & 4023399307\\ 
        \hline
        18 & 324 & 59 & 312418 & 97605006724 & 4023423948\\ 
        \hline
        19 & 361 & 65 & 312419 & 97605631561 & 4023448710\\ 
        \hline
        20 & 400 & 70 & 312420 & 97606256400 & 4023473393\\ 
        \hline
        21 & 441 & 77 & 312421 & 97606881241 & 4023498103\\ 
        \hline
        22 & 484 & 84 & 312422 & 97607506084 & 4023522895\\
        \hline 
    \end{tabular}}
\end{table}

\section{Experiment}
After successful observation of the plots, several equations were experimented with to fit the curves. Details about the experimental equations to fit both the prime count curve in Figure \ref{fig:cnt_prime_vs_num} and difference curve in Figure \ref{fig:diff} are discussed in the following subsections. 

\subsection{Count of primes Curve}
\label{cnt_eqn}
To create quite a similar curve like Figure \ref{fig:cnt_prime_vs_num} a custom equation and four other functions are applied. Also, an equation derived from Bertrand's Postulate is used to evaluate the performance of the new observation compared to the existent ones. This section comprises the explanation for the following functions with their parameters for rivaling the destination curve. 

\begin{enumerate}
    \item Hyperbolic ($cosh()$)
    \item Power Series
    \item Second Degree Polynomial
    \item Conics (From General Equation of Second Degree Curve)
    \item Custom Ratio Equation
    \item Derivation from Bertrand's Postulate
\end{enumerate}

The forecast of the number of primes is accomplished by using only $x$ in the range $[x, x^2]$. The prediction result is expressed as $y$. 

\subsubsection{Hyperbolic}
There are multiple hyperbolic functions such as $tanh()$, $sinh()$, $cosh()$, $cosech()$ etc. Among these functions, the plot of $cosh()$ in a view looks much similar to the curve found. Hence the equation can be written as:
\begin{equation}
\label{eqn:hyper1}
    y = cosh(z)
\end{equation}
The best value for $z$ was found using the following equation:
\begin{equation}
\label{eqn:hyper2}
    z = 1.9029ln(x) - 1.2634
\end{equation}

From Equation \ref{eqn:hyper1} and Equation \ref{eqn:hyper2} we can rewrite as:
\begin{equation}
    y = cosh(1.9023ln(x) - 1.2634)
\end{equation}

\subsubsection{Power Series}
The general equation of a Power Series can be written as:
\begin{equation}
    y = ax^b;
\end{equation}
where $a$ and $b$ are constants. Experimental value found for $a$ and $b$ are:
\begin{center}
    $a = 0.141294556371966$ \\
    $b = 1.90234115616265$
\end{center}

\subsubsection{Second Degree Polynomial}
Second Degree polynomials can be expressed using the following general equation:
\begin{equation}
    y = ax^2 + bx + c
\end{equation}
where $a$, $b$ and $c$ are constants and experimental value found for these constants are:
\begin{center}
    $a = 0.0376$ \\
    $b = 1208.1$ \\
    $c = -3 \times 10^7$
\end{center}
A conditional equation is introduced to lessen the amount of error in the range $[x, x^2]$. As the count of primes is never lower than $x$ in the defined range, the equation can be represented as:

\begin{equation}
    y = 
    \begin{cases}
        x & , \text{if } x > ax^2 + bx + c \\
        ax^2 + bx + c & , \text{otherwise}
    \end{cases}
\end{equation}

\subsubsection{Conics}
Conics are special kind of second degree curves. The conics equation is derived from the general equation of second degree curve. The general equation of a second degree curve can be written as:
\begin{equation}
\label{eqn:poly1}
    Ax^2 + Bxy + Cy^2 + Dx + Ey + F = 0
\end{equation}
Equation \ref{eqn:poly1} can be simplified to find value of $y$ as below:

\begin{equation}
    y = \frac{- (Bx + E) \pm \sqrt{(Bx + E)^2 - 4C(Ax^2 + Dx + F)}}{2C}
\end{equation}
However, the experiment indicates that the expected values can be obtained using only "-" in the equation. Also, the count of primes in the range $[x, x^2]$ is neither negative nor lower than $x$. Hence the equation can be modified as:

\begin{equation}
\label{eqn:conic}
    y = 
    \begin{cases}
        x & , \text{if }x > \frac{- (Bx + E) - \sqrt{(Bx + E)^2 - 4C(Ax^2 + Dx + F)}}{2C} \\
        \frac{- (Bx + E) - \sqrt{(Bx + E)^2 - 4C(Ax^2 + Dx + F)}}{2C} & , \text{otherwise}
    \end{cases}
\end{equation}
Here $A$, $B$, $C$, $D$, $E$ and $F$ are constants. The suitable values found by analyzing the data for the constants are:
\begin{center}
    $A = 3.11199927582249 \times 10^{-09}$ \\
    $B = -9.33817244194697 \times 10^{-15}$ \\
    $C = 3.45730472758733 \times 10^{-21}$ \\
    $D = 0.0000515593800268165$ \\
    $E = -7.63287093993319 \times 10^{-08}$ \\
    $F = -1$ \\
\end{center}
From the values of $A$, $B$, $C$, $D$, $E$ and $F$ we can say that Equation \ref{eqn:conic} is an \textit{Ellipse} if the modification for negative value is removed.

\subsubsection{Custom Ratio Equation}
After plotting, the curve of Figure \ref{fig:ratio} looked like a logarithmic curve. So, we can write:

\begin{equation}
\label{eqn:cus}
    \frac{x^2-x}{y} = \kappa
\end{equation}

Here, $\kappa$ denotes a natural logarithmic equation used for fitting the Figure \ref{fig:ratio}. After several experiments the best value for $\kappa$ was estimated using the following equation:

\begin{equation}
\label{eqn:place}
    \kappa \approx 2.0038ln(x) - 1.0932
\end{equation}

From Equation \ref{eqn:cus} and Equation \ref{eqn:place} we can write:
\begin{equation}
\label{eqn:cus_fin}
    \frac{x^2-x}{y} \approx 2.0038ln(x) - 1.0932
\end{equation}
Considering the approximation as equal for the time being, Equation \ref{eqn:cus_fin} can be modified to find the value of $y$ as below:
\begin{equation}
\label{eqn:cus_calc}
    y = \frac{x^2-x}{2.0038ln(x) - 1.0932}
\end{equation}
It is impracticable to ignore the error generated by considering the approximated value as equal value. That is why the value of $y$ is not a precise match but very close in most cases.

\subsubsection{Derivation of Bertrand's Postulate}
The following equation is derived from Bertrand's Postulate to find the lower limit of number of primes in the range $[x, x^2]$ by using the value of $x$ \cite{reserachgate}.
\begin{equation}
    \frac{1}{2} \log_{2}x^2 = y
\end{equation}
The derived equation claims to find the lower limit always. For evaluation purposes, it is considered that the equation predicts the actual count.

\subsection{Difference Curve}
Although the value of the difference curve fluctuates from time to time, it has a trend to go upward. The goal is to fit a $Straight$ $Line$ for the limited data. The equation used for the experimentation is:
\begin{equation}
    y = 0.0755x + 1018.8
\end{equation}
where $y$ denotes the difference in the number of primes between the ranges $[x, x^2]$ and $[(x-1), (x-1)^2]$.

\section{Result Evaluation}
The performance evaluation of the above-described equations is described in the following subsections. For performance evaluation average relative error is introduced. Both the count of primes and the difference curve bestowed encouraging outcomes.

\subsection{Count of Primes Prediction Evaluation}
\begin{figure}[htbp]
    \centering
    \includegraphics[width=\textwidth,height=\textheight,keepaspectratio]{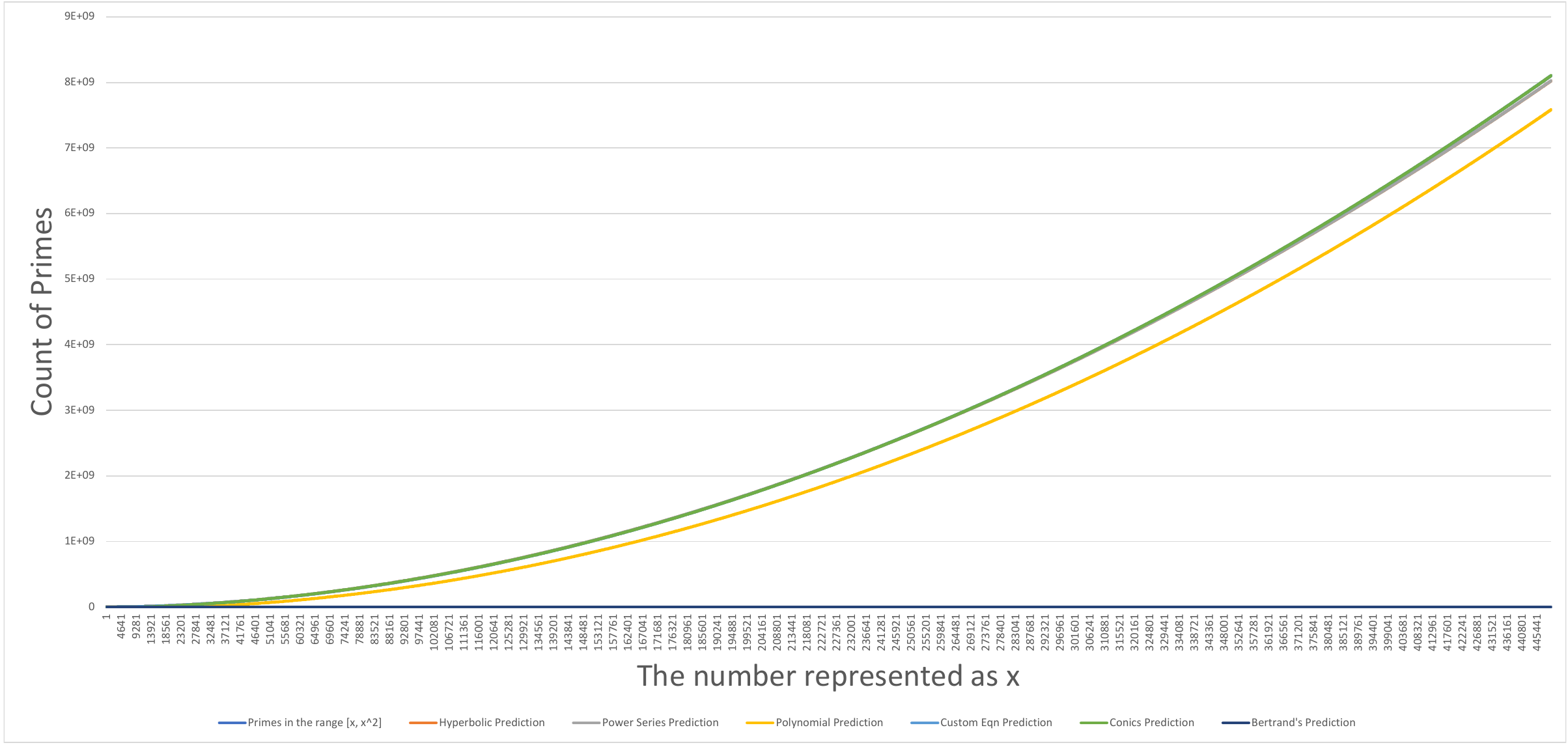}
    \caption{\textbf{Count of Primes Vs. Value of x for Different Functions}}
    \label{fig:comp}
\end{figure}

A comparative chart like Figure \ref{fig:comp} is generated when the predictions by functions are plotted against the real values. The \textbf{Average Relative Error} of all the techniques including Bertrand's Postulate derivation are demonstrated in Table \ref{tab:re}. It is quite imperative to mention that the higher the value of $x$ is, the lower the relative error is in most of the functions.\\

\begin{table}[htbp]
    \centering
    \caption{\textbf{Average Relative Error of all Functions}}
    \label{tab:re}
    \begin{tabular}{|c | c|}
         \hline
         \multicolumn{2}{|c|}{\textbf{Average Value of Relative Errors}}\\ [0.5ex] 
         \hline\hline
         Hyperbolic Equation & 0.64\% \\ \hline
         Power Series Equation & 0.63\% \\ \hline
         Polynomial Equation & 22.31\% \\ \hline
         Custom Ratio Equation & 0.01\% \\ \hline
         Conics Equation & 3.87\% \\ \hline
         Bertrand's Equation & 100.00\% \\ \hline
    \end{tabular}
\end{table}

The values obtained from the equations are very close to actual values. It is also obvious from Figure \ref{fig:comp} and Table \ref{tab:re} that some polishing is still required for accurate values. The statistics of the accurately predicted value is demonstrated in Table \ref{tab:match}. A few correctly predicted values by Custom Equation with floor values are presented in Table \ref{tab:matchExact}.

\begin{table}[htbp]
    \centering
    \caption{\textbf{Predictions using Different Functions}}
    \label{tab:match}
    \resizebox{\linewidth}{!}{\begin{tabular}{| c | c | c | c |}
         \hline
         \textbf{Function} & \textbf{Accurately Predicted} & \textbf{Matched After Ceiling Prediction} & \textbf{Matched After Flooring Prediction}\\ [0.5ex] 
         \hline\hline
         Hyperbolic & 0 & 1 & 0 \\ \hline
         Power Series & 0 & 0 & 2 \\ \hline
         Polynomial & 3 & 3 & 3 \\ \hline
         Custom Ratio & 1 & 28 & 32 \\ \hline
         Conics & 3 & 3 & 3 \\ \hline
    \end{tabular}}
\end{table}

\begin{table}[htbp]
    \centering
    \caption{\textbf{Some Correctly Predicted Values Using Custom Ratio Equation}}
    \label{tab:matchExact}
    \resizebox{\linewidth}{!}{\begin{tabular}{| c | c | c | c |}
         \hline
         \textbf{The Number} & \textbf{Primes in the Range} & \textbf{Prediction from Custom Ratio Equation} & \textbf{Flooring of Predicted Value}\\ [0.5ex] 
         \hline\hline
        731 & 44026 & 44026.3870890 & 44026 \\ \hline
        768 & 48205 & 48205.7312063 & 48205 \\ \hline
        783 & 49949 & 49949.9718114 & 49949 \\ \hline
        858 & 59100 & 59100.2384286 & 59100 \\ \hline
        860 & 59353 & 59353.9894804 & 59353 \\ \hline
        901 & 64666 & 64666.8053193 & 64666 \\ \hline
        922 & 67469 & 67469.6790764 & 67469 \\ \hline
        923 & 67604 & 67604.5254275 & 67604 \\ \hline
        1008 & 79521 & 79521.6327343 & 79521 \\ \hline
        1010 & 79812 & 79812.8299832 & 79812 \\ \hline
        1012 & 80104 & 80104.5217502 & 80104 \\ \hline
        1078 & 90007 & 90007.0211063 & 90007 \\ \hline
        1111 & 95158 & 95158.9030658 & 95158 \\ \hline
        1117 & 96109 & 96109.9237884 & 96109 \\ \hline
        1190 & 108032 & 108032.0843019 & 108032 \\ \hline
        1273 & 122372 & 122372.1181223 & 122372 \\ \hline
        1347 & 135856 & 135856.1727039 & 135856 \\ \hline
        223904 & 2125030424 & 2125030424.8942800 & 2125030424 \\ \hline
        250640 & 2637553503 & 2637553503.8749900 & 2637553503 \\ \hline
        262873 & 2889714221 & 2889714221.3801200 & 2889714221 \\ \hline
        263374 & 2900276799 & 2900276799.6599000 & 2900276799 \\ \hline
        278149 & 3220084404 & 3220084404.7791900 & 3220084404 \\ \hline
        278373 & 3225056401 & 3225056401.8937300 & 3225056401 \\ \hline
        281248 & 3289198019 & 3289198019.5063200 & 3289198019 \\ \hline
        281249 & 3289220435 & 3289220435.1353900 & 3289220435 \\ \hline
    \end{tabular}}
\end{table}

Some predicted values by different equations with relative error against true real values are presented in Table \ref{tab:fval1} and Table \ref{tab:fval2}.

\begin{table}[htbp]
    \centering
    \caption{\textbf{Prime Count using Different Functions with Relative Error against True Values (Part 2}}
    \label{tab:fval1}
    \resizebox{\textwidth}{!}{\begin{tabular}{| c | c | c | c | c | c | c | c |}
         \hline
         \textbf{The Number} & \textbf{Prime Count} & \textbf{Hyperbolic} & \textbf{Hyperbolic R.E.} & \textbf{P.S.} & \textbf{P.S. R.E.} & \textbf{Polynomial} & \textbf{Polynomial R.E.} \\ [0.5ex] 
         \hline\hline
        140001 & 865334106 & 870497682.6 & 0.60\% & 870607669.3 & 0.61\% & 707139663.2 & 18.28\% \\ \hline
        140002 & 865345955 & 870509510.7 & 0.60\% & 870619499.2 & 0.61\% & 707150192.6 & 18.28\% \\ \hline
        140003 & 865357733 & 870521338.9 & 0.60\% & 870631329.2 & 0.61\% & 707160722 & 18.28\% \\ \hline
        140004 & 865369626 & 870533167.3 & 0.60\% & 870643159.2 & 0.61\% & 707171251.4 & 18.28\% \\ \hline
        140005 & 865381494 & 870544995.6 & 0.60\% & 870654989.3 & 0.61\% & 707181781 & 18.28\% \\ \hline
        140006 & 865393246 & 870556824.1 & 0.60\% & 870666819.5 & 0.61\% & 707192310.6 & 18.28\% \\ \hline
        140007 & 865404956 & 870568652.6 & 0.60\% & 870678649.8 & 0.61\% & 707202840.3 & 18.28\% \\ \hline
        140008 & 865416863 & 870580481.2 & 0.60\% & 870690480.2 & 0.61\% & 707213370.1 & 18.28\% \\ \hline
        140009 & 865428765 & 870592309.9 & 0.60\% & 870702310.6 & 0.61\% & 707223899.9 & 18.28\% \\ \hline
        140010 & 865440587 & 870604138.7 & 0.60\% & 870714141.1 & 0.61\% & 707234429.8 & 18.28\% \\ \hline
        140011 & 865452328 & 870615967.5 & 0.60\% & 870725971.7 & 0.61\% & 707244959.8 & 18.28\% \\ \hline
        140012 & 865464130 & 870627796.4 & 0.60\% & 870737802.4 & 0.61\% & 707255489.9 & 18.28\% \\ \hline
        140013 & 865475967 & 870639625.4 & 0.60\% & 870749633.1 & 0.61\% & 707266020.1 & 18.28\% \\ \hline
        140014 & 865487800 & 870651454.5 & 0.60\% & 870761463.9 & 0.61\% & 707276550.3 & 18.28\% \\ \hline
        140015 & 865499530 & 870663283.6 & 0.60\% & 870773294.8 & 0.61\% & 707287080.6 & 18.28\% \\ \hline
        140016 & 865511364 & 870675112.9 & 0.60\% & 870785125.8 & 0.61\% & 707297611 & 18.28\% \\ \hline
        140017 & 865523219 & 870686942.1 & 0.60\% & 870796956.9 & 0.61\% & 707308141.4 & 18.28\% \\ \hline
        140018 & 865534897 & 870698771.5 & 0.60\% & 870808788 & 0.61\% & 707318671.9 & 18.28\% \\ \hline
        140019 & 865546880 & 870710601 & 0.60\% & 870820619.2 & 0.61\% & 707329202.5 & 18.28\% \\ \hline
        140020 & 865558555 & 870722430.5 & 0.60\% & 870832450.5 & 0.61\% & 707339733.2 & 18.28\%  \\ \hline
        140021 & 865570490 & 870734260.1 & 0.60\% & 870844281.8 & 0.61\% & 707350264 & 18.28\% \\ \hline
        140022 & 865582320 & 870746089.8 & 0.60\% & 870856113.2 & 0.61\% & 707360794.8 & 18.28\% \\ \hline
        140023 & 865594219 & 870757919.5 & 0.60\% & 870867944.7 & 0.61\% & 707371325.7 & 18.28\% \\ \hline
        140024 & 865606071 & 870769749.3 & 0.60\% & 870879776.3 & 0.61\% & 707381856.7 & 18.28\% \\ \hline
        140025 & 865617938 & 870781579.3 & 0.60\% & 870891608 & 0.61\% & 707392387.7 & 18.28\% \\ \hline
        140026 & 865629830 & 870793409.2 & 0.60\% & 870903439.7 & 0.61\% & 707402918.8 & 18.28\% \\ \hline
        140027 & 865641590 & 870805239.3 & 0.60\% & 870915271.5 & 0.61\% & 707413450 & 18.28\% \\ \hline
        140028 & 865653402 & 870817069.4 & 0.60\% & 870927103.4 & 0.61\% & 707423981.3 & 18.28\% \\ \hline
        140029 & 865665172 & 870828899.6 & 0.60\% & 870938935.4 & 0.61\% & 707434512.7 & 18.28\% \\ \hline
        140030 & 865677107 & 870840729.9 & 0.60\% & 870950767.4 & 0.61\% & 707445044.1 & 18.28\% \\ \hline
        140031 & 865688926 & 870852560.3 & 0.60\% & 870962599.5 & 0.61\% & 707455575.6 & 18.28\% \\ \hline
        140032 & 865700765 & 870864390.7 & 0.60\% & 870974431.7 & 0.61\% & 707466107.2 & 18.28\% \\ \hline
        140033 & 865712619 & 870876221.2 & 0.60\% & 870986264 & 0.61\% & 707476638.8 & 18.28\% \\ \hline
        140034 & 865724330 & 870888051.8 & 0.60\% & 870998096.3 & 0.61\% & 707487170.5 & 18.28\% \\ \hline
        140035 & 865736181 & 870899882.5 & 0.60\% & 871009928.7 & 0.61\% & 707497702.3 & 18.28\% \\ \hline
        140036 & 865747998 & 870911713.2 & 0.60\% & 871021761.2 & 0.61\% & 707508234.2 & 18.28\% \\ \hline
        140037 & 865759749 & 870923544.1 & 0.60\% & 871033593.8 & 0.61\% & 707518766.2 & 18.28\% \\ \hline
        140038 & 865771569 & 870935375 & 0.60\% & 871045426.4 & 0.61\% & 707529298.2 & 18.28\% \\ \hline
        140039 & 865783391 & 870947205.9 & 0.60\% & 871057259.2 & 0.61\% & 707539830.3 & 18.28\% \\ \hline
        140040 & 865795117 & 870959037 & 0.60\% & 871069092 & 0.61\% & 707550362.5 & 18.28\% \\ \hline
        140041 & 865806830 & 870970868.1 & 0.60\% & 871080924.8 & 0.61\% & 707560894.7 & 18.28\% \\ \hline
        140042 & 865818584 & 870982699.3 & 0.60\% & 871092757.8 & 0.61\% & 707571427.1 & 18.28\% \\ \hline
        140043 & 865830466 & 870994530.6 & 0.60\% & 871104590.8 & 0.61\% & 707581959.5 & 18.28\% \\ \hline
        140044 & 865842310 & 871006361.9 & 0.60\% & 871116423.9 & 0.61\% & 707592492 & 18.28\% \\ \hline
        140045 & 865854121 & 871018193.4 & 0.60\% & 871128257.1 & 0.61\% & 707603024.5 & 18.28\% \\ \hline
        140046 & 865865930 & 871030024.9 & 0.60\% & 871140090.4 & 0.61\% & 707613557.1 & 18.28\% \\ \hline
        140047 & 865877659 & 871041856.5 & 0.60\% & 871151923.7 & 0.61\% & 707624089.8 & 18.28\% \\ \hline
        140048 & 865889475 & 871053688.1 & 0.60\% & 871163757.1 & 0.61\% & 707634622.6 & 18.28\% \\ \hline
        140049 & 865901323 & 871065519.8 & 0.60\% & 871175590.6 & 0.61\% & 707645155.5 & 18.28\% \\ \hline
        140050 & 865913132 & 871077351.7 & 0.60\% & 871187424.1 & 0.61\% & 707655688.4 & 18.28\% \\ \hline
   \end{tabular}}
\end{table}

\begin{table}[htbp]
    \centering
    \caption{\textbf{Prime Count using Different Functions with Relative Error against True Values (Part 1}}
    \label{tab:fval2}
    \resizebox{\textwidth}{!}{\begin{tabular}{| c | c | c | c | c | c | c | c |}
         \hline
        \textbf{The Number} & \textbf{Prime Count} & \textbf{C.R.Eqn} & \textbf{C.R.Eqn R.E.} & \textbf{Conics} & \textbf{Conics R.E.} & \textbf{B.P.D.} & \textbf{B.P.D. R.E.}\\ [0.5ex] 
         \hline\hline
        140001 & 865334106 & 865323992 & 0.00\% & 865796268.5 & 0.05\% & 17.09507761 & 100.00\% \\ \hline
        140002 & 865345955 & 865335807 & 0.00\% & 865808053.1 & 0.05\% & 17.09508791 & 100.00\% \\ \hline
        140003 & 865357733 & 865347622 & 0.00\% & 865819837.8 & 0.05\% & 17.09509822 & 100.00\% \\ \hline
        140004 & 865369626 & 865359437.2 & 0.00\% & 865831622.6 & 0.05\% & 17.09510852 & 100.00\% \\ \hline
        140005 & 865381494 & 865371252.4 & 0.00\% & 865843407.5 & 0.05\% & 17.09511883 & 100.00\% \\ \hline
        140006 & 865393246 & 865383067.7 & 0.00\% & 865855192.5 & 0.05\% & 17.09512913 & 100.00\% \\ \hline
        140007 & 865404956 & 865394883 & 0.00\% & 865866977.5 & 0.05\% & 17.09513943 & 100.00\% \\ \hline
        140008 & 865416863 & 865406698.4 & 0.00\% & 865878762.6 & 0.05\% & 17.09514974 & 100.00\% \\ \hline
        140009 & 865428765 & 865418514 & 0.00\% & 865890547.8 & 0.05\% & 17.09516004 & 100.00\% \\ \hline
        140010 & 865440587 & 865430329.6 & 0.00\% & 865902333.1 & 0.05\% & 17.09517035 & 100.00\% \\ \hline
        140011 & 865452328 & 865442145.2 & 0.00\% & 865914118.5 & 0.05\% & 17.09518065 & 100.00\% \\ \hline
        140012 & 865464130 & 865453961 & 0.00\% & 865925903.9 & 0.05\% & 17.09519096 & 100.00\% \\ \hline
        140013 & 865475967 & 865465776.8 & 0.00\% & 865937689.4 & 0.05\% & 17.09520126 & 100.00\% \\ \hline
        140014 & 865487800 & 865477592.7 & 0.00\% & 865949475 & 0.05\% & 17.09521156 & 100.00\% \\ \hline
        140015 & 865499530 & 865489408.7 & 0.00\% & 865961260.6 & 0.05\% & 17.09522187 & 100.00\% \\ \hline
        140016 & 865511364 & 865501224.7 & 0.00\% & 865973046.4 & 0.05\% & 17.09523217 & 100.00\% \\ \hline
        140017 & 865523219 & 865513040.9 & 0.00\% & 865984832.2 & 0.05\% & 17.09524248 & 100.00\% \\ \hline
        140018 & 865534897 & 865524857.1 & 0.00\% & 865996618.1 & 0.05\% & 17.09525278 & 100.00\% \\ \hline
        140019 & 865546880 & 865536673.4 & 0.00\% & 866008404 & 0.05\% & 17.09526308 & 100.00\% \\ \hline
        140020 & 865558555 & 865548489.7 & 0.00\% & 866020190.1 & 0.05\% & 17.09527339 & 100.00\% \\ \hline
        140021 & 865570490 & 865560306.2 & 0.00\% & 866031976.2 & 0.05\% & 17.09528369 & 100.00\% \\ \hline
        140022 & 865582320 & 865572122.7 & 0.00\% & 866043762.4 & 0.05\% & 17.09529399 & 100.00\% \\ \hline
        140023 & 865594219  & 865583939.3 & 0.00\% & 866055548.7 & 0.05\% & 17.0953043 & 100.00\% \\ \hline
        140024 & 865606071 & 865595756 & 0.00\% & 866067335 & 0.05\% & 17.0953146 & 100.00\% \\ \hline
        140025 & 865617938 & 865607572.7 & 0.00\% & 866079121.5 & 0.05\% & 17.0953249 & 100.00\% \\ \hline
        140026 & 865629830 & 865619389.5 & 0.00\% & 866090908 & 0.05\% & 17.09533521 & 100.00\% \\ \hline
        140027 & 865641590 & 865631206.5 & 0.00\% & 866102694.6 & 0.05\% & 17.09534551 & 100.00\% \\ \hline
        140028 & 865653402 & 865643023.4 & 0.00\% & 866114481.2 & 0.05\% & 17.09535581 & 100.00\% \\ \hline
        140029 & 865665172 & 865654840.5 & 0.00\% & 866126268 & 0.05\% & 17.09536611 & 100.00\% \\ \hline
        140030 & 865677107 & 865666657.6 & 0.00\% & 866138054.8 & 0.05\% & 17.09537642 & 100.00\% \\ \hline
        140031 & 865688926 & 865678474.9 & 0.00\% & 866149841.7 & 0.05\% & 17.09538672 & 100.00\% \\ \hline
        140032 & 865700765 & 865690292.1 & 0.00\% & 866161628.6 & 0.05\% & 17.09539702 & 100.00\% \\ \hline
        140033 & 865712619 & 865702109.5 & 0.00\% & 866173415.7 & 0.05\% & 17.09540733 & 100.00\% \\ \hline
        140034 & 865724330 & 865713927 & 0.00\% & 866185202.8 & 0.05\% & 17.09541763 & 100.00\% \\ \hline
        140035 & 865736181 & 865725744.5 & 0.00\% & 866196990 & 0.05\% & 17.09542793 & 100.00\% \\ \hline
        140036 & 865747998 & 865737562.1 & 0.00\% & 866208777.3 & 0.05\% & 17.09543823 & 100.00\% \\ \hline
        140037 & 865759749 & 865749379.8 & 0.00\% & 866220564.7 & 0.05\% & 17.09544853 & 100.00\% \\ \hline
        140038 & 865771569 & 865761197.5 & 0.00\% & 866232352.1 & 0.05\% & 17.09545884 & 100.00\% \\ \hline
        140039 & 865783391 & 865773015.4 & 0.00\% & 866244139.6 & 0.05\% & 17.09546914 & 100.00\% \\ \hline
        140040 & 865795117 & 865784833.3 & 0.00\% & 866255927.2 & 0.05\% & 17.09547944 & 100.00\% \\ \hline
        140041 & 865806830 & 865796651.3 & 0.00\% & 866267714.9 & 0.05\% & 17.09548974 & 100.00\% \\ \hline
        140042 & 865818584 & 865808469.3 & 0.00\% & 866279502.6 & 0.05\% & 17.09550005 & 100.00\% \\ \hline
        140043 & 865830466 & 865820287.5 & 0.00\% & 866291290.5 & 0.05\% & 17.09551035 & 100.00\% \\ \hline
        140044 & 865842310 & 865832105.7 & 0.00\% & 866303078.4 & 0.05\% & 17.09552065 & 100.00\% \\ \hline
        140045 & 865854121 & 865843924 & 0.00\% & 866314866.3 & 0.05\% & 17.09553095 & 100.00\% \\ \hline
        140046 & 865865930 & 865855742.4 & 0.00\% & 866326654.4 & 0.05\% & 17.09554125 & 100.00\% \\ \hline
        140047 & 865877659 & 865867560.8 & 0.00\% & 866338442.5 & 0.05\% & 17.09555155 & 100.00\% \\ \hline
        140048 & 865889475 & 865879379.3 & 0.00\% & 866350230.7 & 0.05\% & 17.09556186 & 100.00\% \\ \hline
        140049 & 865901323 & 865891197.9 & 0.00\% & 866362019 & 0.05\% & 17.09557216 & 100.00\% \\ \hline
        140050 & 865913132 & 865903016.6 & 0.00\% & 866373807.4 & 0.05\% & 17.09558246 & 100.00\% \\ \hline
    \end{tabular}}
\end{table}

\subsection{Difference Prediction Evaluation}
On average, an \textbf{Average Relative Error} of $12.01\%$ is obtained. The difference also has shown properties that strongly suggest following a distinctive trend. Also, the larger values fitted much well than smaller values.

\section{Conclusion}
The contents of this paper firmly suggest that there is some visible pattern among prime numbers whereas functions like Riemann Prime Counting Function \cite{weisstein_riemann_nodate} exhibit no pattern in the count of prime numbers. The pattern noticeable within a range strongly amplifies the possibility of a pattern in the sequence of prime numbers among themselves. Without such hidden property, it is not perceivable that such a smooth curve will randomly come into existence. If the pattern persists for larger numbers and significant enhancements are performed for accurate prediction it might bring revolution in the sector of prime numbers. Innumerable applications that use prime numbers may get influenced. Many well-established industries including computer-based systems might just be yielded a broad new field of research experimentation. In this paper, the approach was for finding count of primes within a specific range, whereas most prime counting experiments are constructed from $0$ to a specific range. For example, the prime counting function $\pi$ is used to count all the primes within the range from 1 to $x$ in case of $\pi(x)$ whereas the conducted experiment works within the closed range $[x, x^2$. This paper exhibits only some observations and in future more statistical and concrete proofs can be made in this field.

\bibliographystyle{unsrt}  
\bibliography{arXiv} 

\begin{thebibliography}{1}

\bibitem{jjo'connorefrobertson2018}
E~F~Robertson J~J~O'Connor.
\newblock Prime numbers.
\newblock \url{https://mathshistory.st-andrews.ac.uk/HistTopics/Prime_numbers},
  jan 2018.

\bibitem{Lemke}
Robert~J. Lemke~Oliver and Kannan Soundararajan.
\newblock Unexpected biases in the distribution of consecutive primes.
\newblock {\em Proceedings of the National Academy of Sciences},
  113(31):E4446--E4454, 2016.

\bibitem{Klarreich}
Erica Klarreich.
\newblock Mathematicians discover prime conspiracy.

\bibitem{reserachgate}
Issam Kaddoura.
\newblock How many prime numbers between $n^2$ and $(n+1)^2$?, 05 2019.

\bibitem{anaconda}
Anaconda software distribution, 2020.

\bibitem{reback2020pandas}
The pandas~development team.
\newblock pandas-dev/pandas: Pandas, February 2020.

\bibitem{harris2020array}
Charles~R. Harris, K.~Jarrod Millman, St{'{e}}fan~J. van~der Walt, Ralf
  Gommers, Pauli Virtanen, David Cournapeau, Eric Wieser, Julian Taylor,
  Sebastian Berg, Nathaniel~J. Smith, Robert Kern, Matti Picus, Stephan Hoyer,
  Marten~H. van Kerkwijk, Matthew Brett, Allan Haldane, Jaime~Fern{'{a}}ndez
  del R{'{\i}}o, Mark Wiebe, Pearu Peterson, Pierre G{'{e}}rard-Marchant, Kevin
  Sheppard, Tyler Reddy, Warren Weckesser, Hameer Abbasi, Christoph Gohlke, and
  Travis~E. Oliphant.
\newblock Array programming with {NumPy}.
\newblock {\em Nature}, 585(7825):357--362, September 2020.

\bibitem{10.7717/peerj-cs.103}
Aaron Meurer, Christopher~P. Smith, Mateusz Paprocki, Ond\v{r}ej
  \v{C}ert\'{i}k, Sergey~B. Kirpichev, Matthew Rocklin, AMiT Kumar, Sergiu
  Ivanov, Jason~K. Moore, Sartaj Singh, Thilina Rathnayake, Sean Vig, Brian~E.
  Granger, Richard~P. Muller, Francesco Bonazzi, Harsh Gupta, Shivam Vats,
  Fredrik Johansson, Fabian Pedregosa, Matthew~J. Curry, Andy~R. Terrel,
  \v{S}t\v{e}p\'{a}n Rou\v{c}ka, Ashutosh Saboo, Isuru Fernando, Sumith Kulal,
  Robert Cimrman, and Anthony Scopatz.
\newblock Sympy: symbolic computing in python.
\newblock {\em PeerJ Computer Science}, 3:e103, January 2017.

\bibitem{weisstein_riemann_nodate}
Eric~W. Weisstein.
\newblock Riemann {Prime} {Counting} {Function}.

\end{thebibliography}

\end{document}